\documentclass[11pt,a4paper,reqno]{amsart}
\usepackage{amssymb,amsmath,amsfonts,hyperref,amsthm,graphicx,color,float,
mathtools,bbm,caption,subcaption}

\usepackage{tikz}
\usepackage{todonotes}
\usepackage[normalem]{ulem}
\usepackage{fourier}

\numberwithin{figure}{section}
\numberwithin{equation}{section}

\newtheorem{theorem}{Theorem}

\newtheorem{lemma}[theorem]{Lemma}
\newtheorem{proposition}[theorem]{Proposition}

\newtheorem*{conjecture}{Conjecture}

\numberwithin{theorem}{section}

\newcommand{\R}     {\mathbb{R}}
\newcommand{\Rd}     {\mathbb{R}^d}
\newcommand{\N}     {\mathbb{N}}
\renewcommand{\P}   {\mathbb{P}}
\newcommand{\Pp}   {\mathbb{P}_p}

\newcommand{\E}     {\mathbb{E}}
\newcommand{\Ed}     {{\mathbb{E}^d}}
\newcommand{\Ld}     {{\mathbb{L}^d}}
\newcommand{\Z}     {\mathbb{Z}}
\newcommand{\Zd}    {{\mathbb{Z}^d}}

\renewcommand{\d}{{\rm d}}
\newcommand{\indic}[1]	{\mathbbm{1}_{\{#1\}}}
\newcommand{\orig}	{\mathbf{0}}
\newcommand{\conn}{{\leftrightarrow}}
\newcommand{\cluster}{{\mathcal C}}

\newcommand{\IIC}	{\operatorname{IIC}}
\newcommand{\Piic}	{{\mathbb{P}_{\scriptscriptstyle \IIC}}}

\newcommand{\df}	{\operatorname{dim}_f}
\newcommand{\ds}	{\operatorname{dim}_s}
\newcommand{\dm}	{\operatorname{dim}_m}
\newcommand{\dH}	{\operatorname{dim}_H}
\newcommand{\dP}	{\operatorname{dim}_p}

\newcommand{\sss} {\scriptscriptstyle}
\newcommand{\dmin} {d_{\min{}}}

\newcommand{\rad}{{\mathcal R}}  
\newcommand{\cA}{{C_\ast}} 

\title{Fractal dimension of discrete sets and percolation}
\author{Markus Heydenreich}
\address{Mathematisches Institut, Universit\"at M\"unchen, Theresienstr.~39, 80333~M\"unchen, Germany}
\email{m.heydenreich@lmu.de}

\date{11 Dec 2019}

\begin{document}
\begin{abstract}
There are various notions of dimension in fractal geometry to characterise (random and non-random) subsets of $\R^d$. In this expository text, we discuss their analogues for infinite subsets of $\Zd$ and, more generally, for infinite graphs. 
We then apply these notions to critical percolation clusters, where the various dimensions have different values. 
\end{abstract}

\maketitle

\section{What is the dimension of a graph?}\label{intro}
\subsection{Motivation} 
There are various notions of dimension for subsets of $\Rd$, see the classical work of Falconer \cite{Falco14} as well as texts by Fraser and Lehrb\"ack in this volume \cite{Frase19, Lehrb19}. Hausdorff dimension is perhaps the most commonly used, other examples are box dimension and Assouad dimension. Any reasonable notion of dimension yields the same value for strictly self-similar sets, but already for affine self-similar sets these values may differ. 
All these notions depend on microscopic properties of the set, i.e.\ local properties. 

In statistical physics, many interesting models give rise to (random) subsets of the lattice $\Zd$ or even general graphs, and therefore ``dimension'' in this context should describe the macroscopic properties of the set rather than the microscopic ones. 

In this expository text, we shall describe and compare three notions of dimension for graphs: \emph{fractal dimension} and \emph{spectral dimension} can be defined for any (connected and locally finite) graph, while the \emph{mass dimension} requires the graph to be embedded in an ``external'' metric space (for our purpose, we can think of $\Rd$ equipped with the Euclidean norm). 
In the second part, we investigate these notions for (high-dimensional) critical percolation as a prime example of a rich and interesting subset of $\Zd$, and we shall see that the three notions of dimension yield different values. 

It appears that different mathematical communities use different vocabulary, and it is one of our aims to draw the connection between the various concepts involved.

\subsection{Preparatory notions}
We start by recalling basic notions from graph theory. 
Let $G=(V,E)$ be a graph with non-empty vertex set $V$ and edge set $E\subset \binom V 2$ and distinguished vertex $\orig\in V$ (``the root''). 
We interpret $G$ as a metric space with \emph{intrinsic metric} (or `graph metric') 
\begin{equation}
	d_G(x,y)=\inf\big\{n\in \N\colon \exists v_1,\dots,v_n\in V \text{ s.t. } 
		\{x,v_1\},\{v_1,v_2\},\dots,\{v_{n-1},v_n=y\}\in E\big\},
\end{equation}
for the shortest number of edges forming a \emph{path} from $x$ to $y$ (including the case that $d_G(x,y)=\infty$ whenever there is no such path). 

We henceforth assume that the graph is \emph{locally finite}, i.e.\ for all $x\in V$:
\begin{equation}
	\deg_G(x):=\sum_{e\in E}\indic{x\in e}<\infty,
\end{equation}
and \emph{connected}, i.e.\  $d_G(x,y)<\infty$ for all $x.y\in V$. 
For $x\in V$ and $n\in\N_0$, we denote by 
	\[B_x(n):=\big\{y\in V: d_G(x,y)\le n\big\}\] 
the ball w.r.t.\ the intrinsic metric $d_G$, and abbreviate $B(n):=B_\orig(n)$ for the ball of the root. We write $\partial B_x(n):=B_x(n)\setminus B_x(n-1)$ for the inner vertex boundary of $B_x(n)$. 

\subsection{Fractal dimension}
The first notion of dimension is the \emph{fractal dimension} (or ``volume growth dimension'') defined as
\begin{equation}\label{eqDefDf}
	\df(G):=\lim_{n\to\infty} \frac{\log |B(n)|}{\log n}
\end{equation}
whenever the limit exists. More generally, we refer to the upper (resp.\ lower) fractal dimension as $\limsup$ (resp., $\liminf$) of \eqref{eqDefDf}. 
The fractal dimension appears to be a very natural concept, and it characterises the structure of $G$ viewed as a metric space.  
In case of existence of the limit \eqref{eqDefDf}, we can write $|B(n)|=n^{\df(G)+o(1)}$. 

\subsection{Spectral dimension}\label{secSpectralDimension}
A second, completely different approach to dimensionality is given through random walks on the graph $G$. To this end, we define the \emph{(simple) random walk} on the (locally finite) graph $G$ as the (discrete-time) stochastic process with $(S_n)_{n\in\N_0}$ with probability measure $P$ and the property that 
\begin{itemize}
\item $P(S_0=\orig)=1$;
\item For all $n\in\N$ and $x,y\in V$:
\[P(S_{n}=x\mid S_{n-1}=y)=\begin{cases}\frac1{\deg_G(y)}\quad&\text{if $d_G(x,y)=1$},\\0&\text{otherwise.}\end{cases}\]
\end{itemize}
In words, the random walk starts at the root at time $n=0$, and in each time step it moves to one of the neighbouring vertices (chosen independently with equal probability). 

We are now interested in the event that the random walk returns to the origin after a given number $2n$ of steps. Indeed, we may use the decay rate of this probability to define the \emph{spectral dimension} of the graph $G$ as 
\begin{equation}\label{eqDefDs}
	\ds(G):=\lim_{n\to\infty}-2\frac{\log P(S_{2n}=\orig)}{\log n}
\end{equation}
provided that the limit exists. 

Mind that we are interested in returning after an \emph{even} number of steps only, the reason for this is that $P(S_{2n}=\orig)>0$ for all $n$ (e.g. by ``reversing'' the first $n$ steps). However, on bipartite graphs, the random walk can return to the origin {only} after an even number of steps, so that automatically $P(S_n=\orig)=0$ whenever $n$ is odd.

Both notions $\df$ and $\ds$ use the special vertex $\orig$ as `base point'. However, it might be easily observed that $\orig$ is not relevant for the dimension (as long as the graph is connected), and any other vertex of $G$ as base point would lead to the same value of $\df$ and $\ds$. 

The spectral dimension is closely linked to the concept of recurrence and transience of a graph, which we introduce next. To this end, we investigate the probability that the random walk \emph{always} returns to its starting point or not. 
We call the graph \emph{recurrent} if this is the case, i.e., if $P(\exists n\in\N\colon S_n=\orig)=1$. Otherwise, we call the graph \emph{transient}. 
\begin{lemma}\label{lemma-RecTrans}
The graph $G$ is recurrent if $\ds<2$, and it is transient if $\ds>2$. 
\end{lemma}
\proof 
A well-known theorem about random walks (e.g.\ Theorem~5.3.1 in \cite{Durre19}) states that the random walk $(S_n)_n$ is transient if and only if $\sum_{n\in\N}P(S_n=\orig)<\infty$. Thus for $\ds<2$, we have 
\[ \sum_{n\in\N}P(S_n=\orig)\ge \sum_{n\in\N}P(S_{2n}=\orig) = \sum_{n\in\N}n^{-\ds/2+o(1)} =\infty.\]
For an upper bound, we use $P(S_{2n+1}=\orig)\le P(S_{2n}=\orig)$ for all $n\in\N$ (cf.\ Lemma 4.1 in \cite{Barlo17}), and thus
\[ \sum_{n\in\N}P(S_n=\orig)\le 1+ \sum_{n\in\N}2P(S_{2n}=\orig) = 1+ 2\sum_{n\in\N}n^{-\ds/2+o(1)},\]
and this is summable whenever $\ds>2$. 
\qed
\medskip

The ``borderline case'' $\ds=2$ relies on the finer asymptotics of $P(S_{2n}=\orig)$, and thus  the limit \eqref{eqDefDs} is too coarse to give an answer. 

The terminology ``spectral dimension'' suggests a connection with the eigenvalues of the graph Laplacian, see for example Rammal and Toulouse \cite{RammaToulo83} for a discussion in the Physics literature. 
For Brownian motion on a class of compact fractals, Kigami and Lapidus \cite{KigamLapid93} prove a rigorous correspondence between the fractal dimension and the spectrum of the associated Laplacian. 

\subsection{Examples} 
An important example in the present text is the hypercubic lattice $\Ld=(\Zd,\Ed)$ with edge set $\Ed=\big\{\{x,y\}\colon |x-y|=1\big\}$. It is easily observed that $\df(\Ld)=\ds(\Ld)=d$; indeed, this property should hold for any meaningful notion of dimension for discrete sets. 
A Cayley graph is a graph that encodes the abstract structure of a (usually finitely generated) group. The class of Cayley graphs is very rich, and includes the hypercubic lattice, homogeneous trees, and many other graphs. 
Gromov \cite{Gromo81} proved that the limit \eqref{eqDefDf} exists as an integer number for every Cayley graph. Hebisch and Saloff-Coste \cite[Thm.\ 5.1]{HebisSalof93} verified that $\df=\ds$ for Cayley graphs. This equality is true for many other classes of graphs. 

\subsection{The escape time exponent}  
A second notion characterising random walks on graphs is the \emph{escape time exponent} $\beta$, which is defined as 
\begin{equation}\label{eqDefBeta}
	E\big[\inf\{n\in\N\colon S_n\in\partial B(n)\}\big]=n^{\beta+o(1)}.
\end{equation}
Thus $\beta$ describes how long it typically takes to reach the boundary of $n$-balls; by $E[\,\cdot\,]$ we denote expectation w.r.t.\ the random walk measure $P$. For the Euclidean lattice $\Ld$ we have $\beta=2$. If $\beta>2$, we speak of \emph{anomalous diffusion}, which relates to the fact that the random walk moves on average much slower than in Euclidean space: after $n$ steps, the random walk is typically at distance $n^{1/\beta}$ from its starting point. An example for anomalous diffusion is random walk on the Sierpinski gasket, for which Barlow and Perkins \cite{BarloPerki88} proved that $\beta=\log5/\log2$. The exponent $\beta$ is closely linked to $\df$ and $\ds$. Indeed, Barlow and Bass \cite{BarloBass99b} prove that $\beta=2\df/\ds$ for any generalized Sierpinski gasket. 
However, all values of $\beta$ in the interval $[2,\df+1]$ are possible, as pointed out by Barlow \cite{Barlo04b}. 

\subsection{Mass dimension}
The graph notions described above are rather versatile tools for abstract graphs. We shall now consider graphs that are embedded into Euclidean space $\Rd$ (by this we mean that $V\subset \Rd$). For our purpose we can be more restrictive and require that $V\subset \Zd$. 
We denote by 
\begin{equation}\label{eqDefQn}
	Q(n):=[-n,n]^d\cap\Zd
\end{equation}
the ball of radius $n$ with respect to the supremum-metric on $\Zd$. 
The \emph{mass dimension} of a graph $G=(V,E) $ is then defined via 
\begin{equation}\label{eqDefDm}
	\dm(G):=\lim_{n\to\infty} \frac{\log |V\cap Q(n)|}{\log n}.
\end{equation}
Mind the difference between $\ds$ and $\dm$: while the former identifies the growth exponent of balls w.r.t.\ the \emph{intrinsic} (graph) metric, the latter measures balls w.r.t.\ the \emph{extrinsic} (Euclidean) metric. This makes no difference for $G=\Ld$, but we will encounter examples, where this is indeed very different. 
The use of the supremum metric in \eqref{eqDefQn} might appear arbitrary, but since all metrics on $\Zd$ are equivalent, they will all lead to the same value of $\dm$. 

\subsection{Other notions of dimension}
In this exposition we focus on the formerly defined dimensions. However, there are various other notions of dimensions for subsets of $\Zd$ (mostly graph analogues of ``continuum dimensions'' for subsets of $\Rd$). We explain two of these notions, which were introduced by Barlow and Taylor \cite{BarloTaylo89, BarloTaylo92}. 

The first definition is the \emph{discrete Hausdorff dimension} $\dH$, which is defined for subsets of $\Zd$ as follows. 
We say that a set $A\subset \Zd$ is a \emph{finite cube} if there exists $x\in\Zd$ and $r\in\N$ such that $A=Q(r)+\{x\}$ (where $+$ is the Minkowski sum). For a finite set $A\subset \Zd$, we denote by 
\[ \rad(A) = \min\{r:A\subset Q(r)+\{x\}\text{ for some }x\in\Zd\} \]
the radius of $A$ as the radius of a covering cube (and put $R=\infty$ if $|A|=\infty$). 
For $\alpha\ge0$, $A,F\subset \Zd$ and $F\neq\varnothing$, we further let 
\begin{equation*}
	\nu_\alpha(A,F):=\min\left\{ \sum_{i=1}^m\left(\frac{\rad(B_i)}{\rad(F)}\right)^\alpha\colon
		B_1,\dots,B_m\text{ are finite cubes and }A\cap F\subset \bigcup_{i=1}^mB_i\right\}.
\end{equation*}
Let 
\begin{equation}
	m_\alpha(A)=\sum_{n=1}^\infty \nu_\alpha\big(A,Q(2^n)\setminus Q(2^{n-1})\big).
\end{equation} 
Mind that $\alpha\mapsto x^\alpha$ is decreasing for $x\in[0,1]$, and so is $m_\alpha(A)$. We finally define the \emph{discrete Hausdorff dimension} 
\begin{equation}
	\dH(A):=\inf\big\{\alpha\ge0\colon m_\alpha(A)<\infty\big\}. 
\end{equation}
The definition of discrete Hausdorff dimension is clearly modelled by its continuous counterpart. Similarly to (and yet different from) the spectral dimension, this notion is closely related to the recurrence and transience of random walks on $A$: 
\begin{proposition}[{Thm.\ 8.3 in \cite{BarloTaylo92}}]
A set $A\subset \Zd$ is recurrent if $\dH(A)>d-2$, and it is transient if $\dH(A)<d-2$. 
\end{proposition}
Comparison with Lemma \ref{lemma-RecTrans} shows that $\dH$ and $\ds$ often differ. 
What is the behaviour if $\dH(A)=d-2$? If $m_{d-2}(A)<\infty$, then the set is transient as well, but no conclusion is possible when  $m_{d-2}(A)=\infty$, because the $\dH$ is not sensitive enough to decide the matter. 

The second example that we discuss here is the \emph{discrete packing dimension} $\dP$. Its continuous analogue is the packing dimension as defined by Taylor and Tricot \cite{TayloTrico85}, which is the same as Kolmogorov's \emph{metric dimension} and Hawke's \emph{entropy dimension}. 
To this end, we let $A,F\subset \Zd$ as before, and $\varepsilon\in(0,1)$. Then we let 
\begin{equation*}
	\mu_\alpha(A,F,\varepsilon)
	:=	\max\left\{ \sum_{i=1}^m\left(\frac{\rad(B_i)}{\rad(F)}\right)^\alpha\colon
		\genfrac{}{}{0pt}{0}{B_1,\dots,B_m \text{ are finite \emph{pairwise disjoint} cubes}}{\text{centered in }A\cap F \text{ s.t.\ }\rad(B_i)\le\rad(F)^{1-\varepsilon}}\right\},
\end{equation*}
and define the ``packing measure''
\begin{equation}
	p_\alpha(A,\varepsilon)=\sum_{n=1}^\infty \mu_\alpha\big(A,Q(2^n)\setminus Q(2^{n-1},\varepsilon)\big).
\end{equation} 
Then the \emph{discrete packing dimension} is defined as 
\begin{equation}
	\dP(A):=\inf\big\{\alpha\ge0\colon p_\alpha(A,\varepsilon)<\infty\text{ for all }\varepsilon\in(0,1)\big\}. 
\end{equation}
Among the results of Barlow and Taylor \cite[Lemma 3.1]{BarloTaylo92} is the following order of the dimensions: If $A\subset \Zd$, then 
\begin{equation}
	0\le \dH(A)\le \dm(A)\le \dP(A) \le d;
\end{equation}
We return to these notions at the end of this text. 

A different approach to the dimensionality of discrete sets has been proposed recently by Bacelli, Haji-Mirsadeghi, and Khezeli \cite{BacelMirsaKheze18a}.

\section{Percolation}
\subsection{Percolation on $\Ld$}
Percolation theory studies the geometry of certain random subgraphs of $\Ld$. Let $p\in[0,1]$ be a parameter of the model, and make edges in $\Ed$ \emph{occupied} with probability $p$ (independently of each other), and otherwise vacant. More formally, 
we consider the probability space $\Omega=\{0,1\}^\Ed$ equipped with the product topology. 
For a percolation configuration $\omega\in\{0,1\}^\Ed$, an edge $b\in\Ed$ is occupied whenever $\omega(b)=1$, and it is vacant whenever $\omega(b)=0$.
We equip this space with a family of product measures $(\Pp)_{ p\in[0,1]}$ chosen such that $\Pp(b \text{ occupied})=p$ for any $b\in\Ed$ and $p\in[0,1]$. 

\begin{figure}[t]
	\begin{subfigure}[b]{.3\textwidth}
		\includegraphics[width=\textwidth]{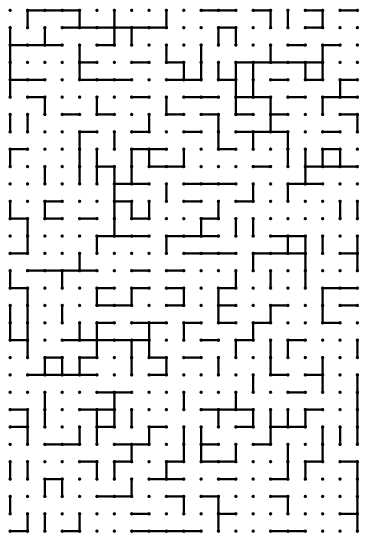}
		\subcaption*{$p=1/3$}
	\end{subfigure}
	\hskip1em
	\begin{subfigure}[b]{.3\textwidth}
		\includegraphics[width=\textwidth]{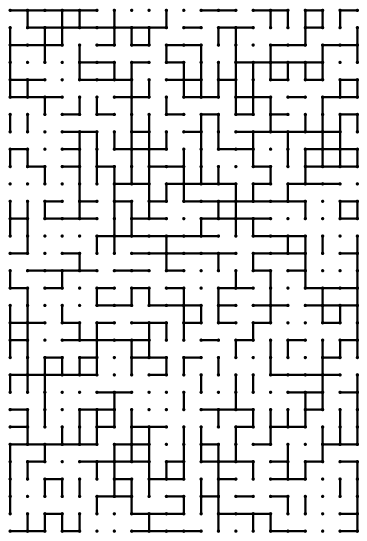}
		\subcaption*{$p=1/2$}
	\end{subfigure}
	\hskip1em
	\begin{subfigure}[b]{.3\textwidth}
		\includegraphics[width=\textwidth]{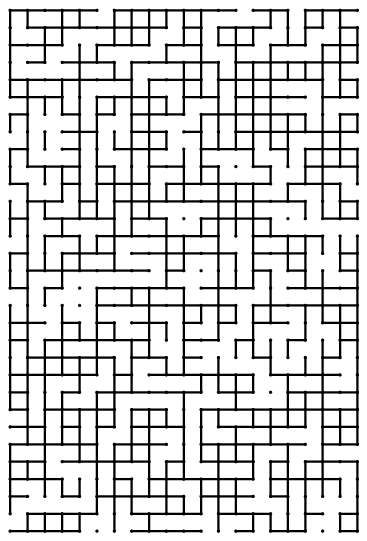}
		\subcaption*{$p=2/3$}
	\end{subfigure}
	\caption{Three realisations of percolation on $\mathbb L^2$.}
	\label{fig-percolation}
\end{figure}

We say that $x$ is {\em connected} to $y$ and write $x\conn y$ when there exists a (finite) path of occupied edges connecting $x$ and $y$. 
Formally,  $x\conn y$ on a configuration $\omega\in\{0,1\}^\Ed$ if there exist $x=v_0,v_1,\dots,v_{m-1},v_m=y\in\Zd$ with the property that $\{v_{i-1},v_i\}\in\Ed$ and $\omega(\{v_{i-1},v_i\})=1$ for all $i=1,\dots,m$ ($m\in\N$). 
We further write $\{x\conn y\}=\{\omega\colon x\conn y\text{ on the configuration }\omega\}$. 
We let the {\em cluster} of $x$ be all the vertices that are connected to $x$, i.e., $\cluster(x)=\{y\colon x\conn y\}$. By convention, $x\in \cluster(x)$.

We define the \emph{percolation function} $p\mapsto \theta(p)$ by
\begin{equation}
    \label{theta-def}
    \theta(p)=\Pp(|\cluster(x)|=\infty),
\end{equation}
where $x\in \Z^d$ is an arbitrary vertex and $|\cluster(x)|$ denotes the number of vertices in $\cluster(x)$.
By translation invariance, the above probability does not depend on the choice of $x$. We therefore often investigate $\cluster = \cluster(\orig)$ where $\orig\in \Z^d$ denotes the origin. 

When $\theta(p)=0$, then the probability that the origin is inside an infinite connected component is 0, so that there is almost surely no infinite connected component. 
On the other hand, when $\theta(p)>0$, then (by ergodicity) the proportion of vertices in infinite connected components equals $\theta(p)>0$, and we say that the system \emph{percolates}. 

We define the \emph{percolation critical value} by
\begin{equation}
    \label{pc-def}
    p_c=\inf\{p\colon \theta(p)>0\}.
\end{equation}
It is well-known that $p_c=1$ on the one-dimensional lattice $\mathbb L^1$ and $p_c\in(0,1)$ on $\Ld$ for all $d\ge2$. For this and other basic properties we refer to the textbooks by Grimmett \cite{Grimm99}, Bollobas and Riordan \cite{BolloRiord06} and Werner \cite{Werne09}. 
See Figure \ref{fig-percolation} for a simulation of percolation with different values of $p$.

For every percolation realization $\omega\in\Omega$, we can define a random walk on the cluster $\cluster$ as in Section \ref{secSpectralDimension}; we denote the corresponding measure by $P^\omega$. Random walk on percolation clusters is a benchmark model of \emph{random walk in (non-elliptic) random environment}. 

\subsection{Dimension of percolation clusters}
We now address the question: What is the dimension of the percolation cluster $\cluster=\cluster(\orig)$? The answer depends on the value of $p$. 

Indeed, if $p<p_c$, then $|\cluster|<\infty$ for $\Pp$-almost all $\omega$, and hence 
\begin{equation}
	\df(\cluster)=\dm(\cluster)=\ds(\cluster)=0\qquad \Pp-a.s.
\end{equation}

We get a different picture when $p>p_c$, and thus $\theta(p)=\Pp(|\cluster(x)|=\infty)>0$. We condition on the event that the origin lies in an infinite cluster, and denote the conditional probability by $\Pp^\ast(\;\cdot\;)=\Pp(\;\cdot\;\mid 0\conn\infty)$. 
It may be seen by applying the ergodic theorem that \begin{equation}
	\lim_{n\to\infty} \frac{|\cluster\cap Q(n)|}{|Q(n)|}\to\theta(p)\qquad \Pp^\ast-\text{a.s.,}
\end{equation}
and hence $\dm(\cluster)=d$. Furthermore, we get that $\df(\cluster)=d$ (almost surely w.r.t.\ the measure $\Pp^\ast$) by exploiting the large deviation bounds on the graphical distance by Antal and Pisztora \cite{AntalPiszt96}. Concerning the spectral dimension, Barlow \cite{Barlo04} proved the heat kernel bounds 
\begin{equation}
	c_1 n^{-d/2}\le P^\omega(S_{2n}=\orig)\le c_2 n^{-d/2}
\end{equation}
for $\Pp^\ast$-almost all $\omega$ (where the constants $c_1,c_2>0$ depend on the value of $p$), and hence $\ds(\cluster)=d$ as well. 
Barlow's result was strengthened further to get a \emph{quenched invariance principle} \cite{BergeBisku07, MathiPiatn07}. 

Finally, the critical case $p=p_c$. 
There is a rather general lower bound on the volume growth of critical clusters. Recall that we denote by $B(n)$ the ball w.r.t.\ the intrinsic (graph) metric $d_\cluster$ on the cluster $\cluster=\cluster(\orig)$, and $\partial B(n)=B(n)\setminus B(n-1)$.  
\begin{theorem}\label{thmLowerBdB}
For percolation on $\Ld$, $d\ge1$, 
	\begin{equation}
	\E_{p_c}|B(n)|\ge n, \qquad n\ge1. \label{eqLowerBdB}
	\end{equation}
\end{theorem}
We provide a proof at the end of this chapter. Mind that \eqref{eqLowerBdB} implies that $\liminf_{n\to\infty}\log \E_{p_c}|B(n)|/\log n\ge1$, and it might be tempting to conjecture that even $\df\ge1$. This, however, is not true. Even stronger, it is strongly believed that critical infinite clusters do not exist: 
\begin{conjecture}
For percolation on $\Zd$, $d\ge2$, we have that $\theta(p_c)=0$, and thus $|\cluster(x)|<\infty$ for all $x\in\Zd$ $\P_{p_c}$-a.s.
\end{conjecture}
The conjecture is known to be true for $d=2$ by Kesten \cite{Keste80} as well as in high dimensions by Hara and Slade \cite{HaraSlade90a}, where the meaning of  \emph{high dimensions}  is that there exists $\dmin{}>6$ such that the claim is true for $d\ge \dmin$. Fitzner and van der Hofstad \cite{FitznHofst17} optimized the strategy of Hara and Slade and verified that $\dmin=11$ suffices. 
Proving this conjecture in dimensions $3\le d\le 10$ is a major open problem in percolation theory; see also \cite{Grimm99} and \cite[Open Problem 1.1]{HeydeHofst17}. 

In view of the presumed result that $\theta(p_c)=0$, we thus get that all clusters are almost surely finite and hence all dimensions equal 0, precisely as for $p<p_c$. Yet an interesting structure will emerge if we look at the interesting geometry of critical clusters from a different angle. We now  investigate this further for the two regimes that we do understand rigorously, namely $d=2$ and \emph{high dimensions}.

\subsection{The incipient infinite cluster}
When $\theta(p_c)=0$, this leaves us with a most remarkable situation:
At the critical point $p_c$ there are clusters at all length scales, which are, however, all finite. As we then make a density $\varepsilon>0$ of closed edges open, the large clusters connect up to form a (unique) infinite cluster, no matter how small $\varepsilon$ is.  At criticality, the critical cluster is therefore \emph{at the verge of appearing}. This observation motivated the introduction of an \emph{incipient infinite cluster} (IIC) as a critical cluster that is \emph{conditioned} to be infinite. 

Somewhat simplified, the \emph{incipient infinite cluster} (IIC) is defined as the cluster of the origin under the critical measure $\P_{p_c}$ conditioned on $\{|\cluster(0)|=\infty\}$. Since this would condition on an event of zero probability, a rigorous construction of the IIC requires a limiting argument. The first mathematical construction has been carried out by Kesten \cite{Keste86a} in two dimensions, who considered two limiting schemes: 
\begin{itemize}
\item[$\rhd$] under $\P_{p_c}$, condition on the event $\{\cluster(0)\cap \partial\Lambda_n\neq\varnothing\}$, and then let $n\to\infty$;
\item[$\rhd$] under $\Pp$ ($p>p_c$), condition on the event $\{|\cluster(0)|=\infty\}$ and let $p\searrow p_c$.
\end{itemize}
Kesten proved that both limits exist in dimension $d=2$, and lead to the {\em same} limiting measure, which he calls the incipient infinite cluster.  
He was motivated by observations in the physics literature, which indicated anomalous diffusion for random walk on large critical percolation clusters. 
Kesten \cite{Keste86b} confirmed this, and proved that the exit time exponent $\beta$ satisfies $\beta>2$ on incipient infinite cluster in two dimensions. It is an open problem to improve this bound. 

For percolation on a regular tree, the cluster distribution is precisely that of a Galton-Watson tree with binomial offspring distribution. 
Hence, the incipient infinite cluster for percolation on a tree is a special case of critical Galton-Watson tree conditioned on non-extinction. 
It was again Kesten \cite{Keste86b} who studied the latter, and proved that it can be constructed in two steps: a single infinite line of descent, casually phrased as ``the immortal particle'' and more formally as ``cluster backbone'', and critical trees hanging off this backbone. 
He further investigated the escape time exponent for this incipient infinite cluster on trees, and proved that $\beta=3$.

We now come to the case of high-dimensional percolation, where the IIC was constructed by van der Hofstad and J\'arai: 
\begin{theorem}[IIC construction \cite{HofstJarai04}]
\label{thm-IIC}
There is a dimension $\dmin>6$ such that for $d\geq \d_{\min}$ and any event $E$ that depends on the status of finitely many edges, the limit 
	\begin{equation}\label{eqPiic}
	\P_{\rm\sss IIC}(E):=\lim_{|x|\to\infty}\P_{p_c}\big(E\mid 0\conn x\big)
	\end{equation}
exists. 
\end{theorem}
The limitation to events that depend on the status of only finitely many edges is a technical one. In fact, such events form an algebra on $\Omega$ which is stable under intersections, and we may thus extend $\Piic$ to a measure on the $\sigma$-fields generated by the product topology. We denote this measure $\Piic$ the \emph{incipient infinite cluster measure}. 

It is straightforward to see that indeed $\P_{\rm\sss IIC}(|\cluster(\orig)|= \infty)=1$, as desired. Since $\theta(p_c)=0$, the IIC is also {\em one-ended} in the sense that the removal of any finite region of the IIC leaves one infinite part. It can be seen that the infinite path is {\em essentially unique} in the sense that any pair of infinite self-avoiding paths in the IIC share infinitely many edges.

Van der Hofstad and J\'arai derive also another construction of the IIC-measure in high dimensions, namely
	\begin{equation}
	\label{eqQiic} 
	\displaystyle
	\mathbb P_{\rm\sss IIC}(E)
	=\lim_{p\nearrow p_c}\frac{\sum_{x\in\Zd}
	\Pp(E\cap\{0\conn x\})}{\sum_{x\in\Zd}\Pp(0\conn x)}. 
	\end{equation}
A third construction (same as Kesten's first construction in two dimension) was derived with van der Hofstad and Hulshof \cite{HeydeHofstHulsh14a}. 

Mind that the measure $\P_{\sss\rm IIC}$ has lost the translation invariance of the percolation measures $\Pp$. Indeed, the origin $\orig$ plays a special role, since we have enforced that the cluster $\cluster(\orig)$ is infinite. 

\subsection{Lower bound for the expected size of critical balls}
We now prove Theorem \ref{thmLowerBdB}. One ingredient is an alternative characterization of $p_c$, namely 
\begin{equation}\label{eqPcAlt}
	p_c=\sup\big\{p\in[0,1]\colon \E_{p}|\cluster|<\infty\big\},
\end{equation}
which is standard in percolation theory \cite{AizenBarsk87,Mensh86}; we also refer to the  short proof by Duminil-Copin and Tassion \cite{DuminTassi16}. 
In our proof of Theorem \ref{thmLowerBdB}, we adapt ideas of \cite{DuminTassi16} but use the intrinsic (graph) metric rather than the extrinsic one. It appears that the proof is valid in much wider context, namely all transitive connected graphs whose percolation threshold is strictly between 0 and 1 and for which \eqref{eqPcAlt} is true. 

\proof[Proof of Theorem \ref{thmLowerBdB}]
We define the value 
$ \bar p_c = \sup M$ with 
\[M=\big\{p\in(0,1)\colon \exists n\in\N \text{ such that }\E_p|\partial B(n)|<1\big\}.\]
Fix an arbitrary $p<\bar p_c$. Then exists $\varepsilon>0$ and $n\in\N$ such that $\E_p|\partial B(n)|<1-\varepsilon$. Fix such $\varepsilon$ and $n$. We now claim that 
\begin{equation}\label{eqBknBd}
	\P_{p}\big(\partial B(kn)\neq 0\big)\le (1-\varepsilon)^k,\qquad k\in\N.
\end{equation}
The proof of \eqref{eqBknBd} is via induction in $k$. The initialization of the induction is our assumption. For the inductive step, we assume that \eqref{eqBknBd} is true for some $k$, and aim to prove it for $k+1$. 
We first condition on the ball $B(n)$: 
\begin{equation}
	\P_{p}\big(\partial B((k+1)n)\neq 0\big)
	=\sum_{A\subset\Zd} \P_{p}\big(B(n)=A,\partial B((k+1)n)\neq 0\big).
\end{equation}
We treat the set $A$ as a subgraph of $\Ld$, and denote by $\partial A$ the vertices with maximal graphical distance from $\orig$, this allows us to bound 
\begin{align}
	\P_{p}\big(\partial B((k+1)n)\neq 0\big)
	&=\sum_{A\subset\Zd} \P_{p}\big(B(n)=A,\bigcup_{y\in\partial A}\partial B_y(kn)\neq 0\text{ in }(\Zd\setminus A)\cup\{y\}\big)\nonumber\\
	&\le \sum_{A\subset\Zd}\sum_{y\in\partial A} \P_{p}\big(B(n)=A,\partial B_y(kn)\neq 0\text{ in }(\Zd\setminus A)\cup\{y\}\big).
\end{align}
The event $\{B(n)=A\}$ depends on the status of all the edges with at least one endpoint in $A\setminus\partial A$. On the other hand, $\big\{\partial B_y(kn)\neq 0\text{ in }\Zd\setminus A\cup\{y\}\big\}$ depends on the status of the edges not touching $A\setminus \partial A$. Hence, the two events are independent, and we bound further 
\begin{align}
	\P_{p}\big(\partial B((k+1)n)\neq 0\big)
	&\le \sum_{A\subset\Zd}\sum_{y\in\partial A} \P_{p}\big(B(n)=A\big)\,\P_{p}\big(\partial B_y(kn)\neq 0\text{ in }(\Zd\setminus A)\cup\{y\}\big)
	\nonumber\\
	&\le \sum_{A\subset\Zd}\sum_{y\in\partial A} \P_{p}\big(B(n)=A\big)\,\P_{p}\big(\partial B_y(kn)\neq 0\big)
\end{align}
Transitivity of the underlying lattice gives $\P_{p}\big(\partial B_y(kn)\neq 0\big)=\P_{p}\big(\partial B(kn)\neq 0\big)$. 
Since 
\[\sum_{A\subset\Zd}\sum_{y\in\partial A} \P_{p}\big(B(n)=A\big)
	=\sum_{A\subset\Zd}|\partial A|\; \P_{p}\big(B(n)=A\big)
	=\E_p|\partial B(n)|
	\le 1-\varepsilon,
\]
we can use the induction hypotheses to obtain $\P_{p}\big(\partial B((k+1)n)\neq 0\big)\le (1-\varepsilon)^{k+1}$, thus proving \eqref{eqBknBd}. 
Consequently, 
\[\P_p(|\cluster|=\infty)\le \lim_{k\to\infty}\P_{p}\big(\partial B(kn)\neq 0\big)=0\]
and thus $p\le p_c$. Since $p<\bar p_c$ was arbitrary, we conclude $\bar p_c\le p_c$. 

We further observe that $M$ is an open subset of $[0,1]$, and therefore $\bar p_c\not\in M$. 
This implies $\E_p|\partial B(n)|\ge1$ for all $n\in\N$, and thus 
\[ \E_{\bar p_c}|\cluster|=\sum_{n\in\N_0}\E_{\bar p_c}|\partial B(n)|\ge\sum_{n\in\N_0}1=\infty,\]
and via \eqref{eqPcAlt} we thus get that $\bar p_c\ge p_c$. 
Together with the foregoing, we established $\bar p_c=p_c$. 

The finishing touch is provided by the partial summation 
\[ \E_{p_c}|B(n)|=\E_{\bar p_c}|B(n)|=\sum_{k=0}^n\E_{\bar p_c}|\partial B(k)|\ge\sum_{k=1}^n1=n.\]
\qed

\section{Dimension of the incipient infinite cluster}

In this section we come to the main endeavour of this text, which is characterising the various dimensions of incipient infinite cluster. 

Let us deal with the planar case first. Kesten \cite{Keste86a} proved for various two-dimensional lattices that 
\begin{equation}
	\lim_{\lambda\to\infty}\Piic\left(\lambda^{-1}\le \frac{|\cluster\cap Q(n)|}{n^2\;\P_{p_c}(\orig\conn\partial Q(n))}\le \lambda	\right)=1
\end{equation} 
uniformly in $n$. 
For the case of site percolation on the triangular lattice, it is known that $P_{p_c}(\orig\conn\partial Q(n))=n^{-5/48+o(1)}$, cf.\ \cite{LawleSchraWerne02}. This suggests that the mass dimension equals $\dm=91/96$. However, in view of Lemma \ref{lemCriterion} below, the control of the error terms is not strong enough to conclude that $\dm=91/96$ in an ($\Piic$-)almost sure sense. 
Concerning the fractal dimension $\df$, it is a challenging open problem to derive sharp bounds on the \emph{intrinsic} (graph) distance of critical two-dimensional clusters. For a recent survey of bounds on the intrinsic distance in the planar case, we refer to Damron \cite{Damro16}.

We now come to the case of high-dimensions, where the results are most complete. For a general survey of results in high-dimensional percolation, we refer to our recent textbook \cite{HeydeHofst17}.

%

For the incipient infinite cluster in high dimensions, the results are summarized in the following theorem: 
\begin{theorem}[{\cite{Cames15,KozmaNachm09}}]\label{thmDimension}
	For the incipient infinite cluster in high dimensions, we have that 
	\[ \ds(\cluster)=4/3, \quad \df(\cluster)=2, \quad \dm(\cluster)=4 \qquad  \Piic-\text{a.s.} \]
\end{theorem}
Interestingly, all three dimensions of $\cluster$ are independent of the dimension $d$ of the embedding space, an indication that the geometry of the embedding space is less visible, and the model appears similar as their non-spatial analogues (as predicted by mean-field theory). 

In the sequel, we demonstrate the proof for the fractal dimension based on a number of standard results for high-dimensional percolation.  
Finally, we discuss the necessary adaptations for the other dimensions $\ds$ and $\dm$. 

The analysis of percolation in high dimension is rooted in a technique called the \emph{lace expansion}. For percolation, this was pioneered in a seminal 1990 paper by Hara and Slade \cite{HaraSlade90a}, who were inspired by earlier work of Brydges and Spencer \cite{BrydgSpenc85} for self-avoiding walk. For our purpose we need the following estimate on the percolation connectivity: there exist $C,c>0$ such that for all $x,y\in\Zd$, $x\neq y$, 
\begin{equation}\label{eqtpf}
	c|x-y|^{d-2}\le \P_{p_c}(x\conn y)\le C|x-y|^{d-2}. 
\end{equation}
The estimate \eqref{eqtpf} was first derived for a spread-out version of percolation \cite{HaraHofstSlade03}, and adapted by Hara \cite{Hara08} to our setting. Fitzner and van der Hofstad \cite{FitznHofst17} verified that it is valid in dimension $d>10$. 

The upper bound in \eqref{eqtpf} readily implies the famous \emph{triangle conditon}, which in turn implies that various \emph{critical exponents} take on their mean-field values. 
We need only two implications here, and refer to a general discussion of critical exponents to \cite[Section 1.2]{HeydeHofst17}: There are constants $c_1,C_1,c_2,C_2>0$ such that 
\begin{equation}\label{eqGamma}
	\frac{c_1}{p_c-p} \le \E_p|\cluster| \le \frac{C_1}{p_c-p},
	\qquad p\in(0,p_c),
\end{equation}
and 
\begin{equation}\label{eqDelta}
	\frac{c_2}{\sqrt{k}} \le \P_{p_c}\big(|\cluster|\ge k) \le \frac{C_2}{\sqrt{k}},
	\qquad k\in\N.
\end{equation}
The bound \eqref{eqGamma} is due to Aizenman and Newman \cite{AizenNewma84}, the bound \eqref{eqDelta} due to Barsky and Aizenman \cite{BarskAizen91}. 

Our final ingredient is the famous BK-inequality. To this end, we define the \emph{disjoint occurrence} $E\circ F$ of two events $E$ and $F$ as 
\begin{equation}\label{eqCirc}
	E\circ F = \big\{ \omega\colon \exists K\subset \Ed\text{ such that }\omega_K\in E, \omega_{\Ed\setminus K}\in F\big\},
\end{equation}
where $\omega_K := \{\omega'\colon \omega(e)=\omega'(e)\text{ for all }e\in K\}$ is the ``$K$-cylinder of $\omega$''. 
Then the BK-inequality \cite{BergKeste85,Reime00} establishes that 
\begin{equation}\label{eqBK}
	\P_p(E\circ F) \le \P_p(E) \;\P_p(F)
\end{equation}
for any $p\in[0,1]$ and all events $E,F$ that \emph{depend on finitely many edges}. This last confinement can be lifted in many cases, and indeed \eqref{eqBK} is true for all events that we are considering in the present text (see also Section 2.3 in \cite{Grimm99}). 

\subsection{The fractal dimension}
We now prove that $\df(\cluster)=2$ whenever \eqref{eqtpf} is valid. 
We start by showing that the lower bound \eqref{eqLowerBdB} has a matching upper bound in high dimensions (which is supposedly false in dimension $d<6$). 
\begin{lemma}[Ball growth]\label{lemBalls}
Consider percolation in dimension $d>10$. 
There exists a constant $C_3>0$ such that for all $n\in\N$, 
\[ n\le \E_{p_c}|B(n)| \le C_3 n.\]
\end{lemma}

\proof
The lower bound was already contained in \eqref{eqLowerBdB}. 
We follow Sapozhnikov \cite{Sapoz10} for a proof of the upper bound. 
Let $p<p_c$. We consider the following coupling of percolation with parameters $p$ and $p_c$: Starting with a critical percolation configuration (edges are occupied with probability $p_c$), make every occupied edge vacant with probability $1-(p/p_c)$. 
This construction implies that for any $x\in\Zd$, $p<p_c$, and $n\in\N$,
	$$  
	\P_{p}\big( d_{\cluster}(\orig,x)\leq n\big)
	\ge \left(\frac {p}{p_c}\right)^n\P_{p_c}(d_{\cluster}(\orig,x)\leq n).
	$$
Summing over $x$ and using the inequality $ \P_{p}(d_{\cluster}(0,x)\leq n)\le \P_{p}(0\conn x)$, we obtain 
$$ \E_{p_c}|B(n)|
	\le\left(\frac {p_c}p\right)^n \E_p|\cluster|
	\le C_1 \left(\frac {p_c}p\right)^n (p_c-p)^{-1},$$
where the last bound comes from \eqref{eqGamma}. 
Choosing $p=p_c(1-\frac1{2n})$ proves the claim. 
%
\qed

\begin{lemma}[{Arm exponents \cite{KozmaNachm09}}]\label{lemArmExponent}
Consider percolation in dimension $d>10$. 
There exists constants $C,c>0$ such that for all $n\in\N$, 
\[ \frac cn\le \P_{p_c}\big(\partial B(n)\neq \varnothing\big) \le \frac {C}n.\]
\end{lemma}
\proof
We start with the proof of the lower bound, and use the well-known second-moment method. 
The basic inequality is 
\begin{equation}\label{eq SecondMoment}
	\P(Z>0)\ge (\E Z)^2/\E Z^2,
\end{equation}
which is valid for any non-negative random variable $Z$. 
We aim to apply this to $Z=\big| B(\lambda n)\setminus B(n)\big|$ with $\lambda=2C_3$. 
Now Lemma \ref{lemBalls} yields 
	\[
	\E_{p_c}|B(\lambda n)\setminus B(n)|
	\ge \lambda n-C_3n=C_3n.	
	\]
We now estimate the second moment of $B(\lambda  n)$. 
Indeed, if both $x$ and $y$ are connected with distance $\le\lambda n$ from $\orig$, then there must exist a ``branch point'' $z\in\Zd$ such that there are (edge-)disjoint paths from $\orig$ to $z$, from $z$ to $x$ and from $z$ to $y$. We may use the symbol $\circ$ (recall \eqref{eqCirc}) to write this as 
\begin{align*}
	&\{d_{\cluster}(\orig,x)\le {\lambda } n\} \cap \{d_{\cluster}(\orig,y)\le {\lambda } n\}  \\
	&\qquad\subseteq \bigcup_{z}\{d_{\cluster}(\orig,z)\le {\lambda } n\} \circ \{d_{\cluster}(z,x)\le {\lambda } n\} \circ \{d_{\cluster}(z,y)\le {\lambda } n\}. 
\end{align*}
Consequently, the BK-inequality \eqref{eqBK} and Lemma \ref{lemBalls} yields 
	\begin{align}
 	\E_{p_c}|B({\lambda } n)|^2
	&= \sum_{x,y} \P_{p_c}(d_{\cluster}(\orig,x)\le {\lambda } n,d_{\cluster}(\orig,y)\le {\lambda } n)\nonumber\\
	&\le \sum_{x,y,z} \P_{p_c}(d_{\cluster}(\orig,z)\le {\lambda } n)\,\P_{p_c}(d_{\cluster}(z,x)\le {\lambda } n)\,\P_{p_c}(d_{\cluster}(z,y)\le {\lambda } n)\nonumber\\
	&= \Big[\sum_{z\in\Zd}\P_{p_c}(d_{\cluster}(\orig,z)\le {\lambda } n)\Big]^3
	=B(\lambda n)^3 
	\le C'n^3,
	\label{eq tree graph bound}
\end{align}
for some constant $C'>0$. 
Consequently, the bound in \eqref{eq SecondMoment} yields 
	\[	
	\P_{p_c}\big(\exists x\in \Z^d\colon d_{\cluster}(0,x)\geq {\lambda } n\big)
	\ge \P_{p_c}\big(|B({\lambda } n)\setminus B(n)|>0\big)
	\ge \frac{C_3^2n^2}{C'n^3}=\frac{C_3^2}{C'n}, 
	\]
which proves the statement with $c=\frac{C_3^2}{\lambda C'}=\frac{C_3}{2C'}$. 

The upper bound uses a clever induction argument. 
For subgraphs $G$ of the infinite lattice $\Ld$, we denote by $\cluster_G=\cluster_G(\orig)$  the (restricted) percolation cluster of $\orig$ \emph{in the subgraph $G$}, and denote by $B_{\cluster_G}(n)=\big\{y\in\Zd\colon d_{\cluster_G}(\orig,y)\le n\big\}$ the corresponding ball w.r.t.\ the graph metric on the restricted cluster $\cluster_G$. 
We further define 
	\[
	H(n;G):=\big\{\partial B_{\cluster_G}(n)\not = \varnothing\big\}  
	\]
for the ``one-arm event'' on the graph $G$, and 
	\[  
	\Gamma(n)= \sup\big\{\P_{p_c}(H(n;G))\colon G\text{ is subgraph of } \Ed\big\}. 
	\]
It turns out that working with $\Gamma(n)$ rather than $\P_{p_c}(H(n;\Ld))$ enables us to apply a regeneration argument, which would not work for $\P_{p_c}(H(n;\Ld))$, since it is not monotone. 

For $C_2$ as in \eqref{eqDelta}, we choose $\cA\ge1$ large enough so that
	\begin{equation}
	\label{eqArmBd2}
	3^3\cA^{2/3}+C_2 \cA^{2/3}\le \cA,
	\end{equation}
We claim that, for any integer $k\ge 0$, 
	\begin{equation}
	\label{eq Gamma A}
	\Gamma(3^k)\le \frac{\cA}{3^k}.
	\end{equation}
This readily implies the upper bound of the lemma, since for any $n$ we choose $k$ such that $3^{k-1}\le n< 3^k$ and then
	\[ 
	\P_{p_c}(H(n;\Ld)) \le \Gamma(n) \le \Gamma(3^{k-1}) \le \frac{\cA}{3^{k-1}} \le \frac{3\cA}{n}.  
	\]
The proof of \eqref{eq Gamma A} is via induction in $k$. The claim is trivial for $k=0$ since $\cA\ge1$. For the inductive step we assume \eqref{eq Gamma A} for $k-1$ and prove it for $k$.  Depending on the size $|\cluster_G|$ of the restricted cluster $\cluster_G$ for arbitrary subgraphs $G$, we estimate
	\begin{equation}
	\label{eqArmBd3}
	\P_{p_c}(H(3^k;G))
	\le \P_{p_c}\big(H(3^k;G), |\cluster_G|\le \cA^{-4/3}9^k\big)
	+ \P_{p_c}\big(|\cluster_G|> \cA^{-4/3}9^k\big).
	\end{equation}
For the second summand, we use \eqref{eqDelta} to obtain
	\begin{equation}	
	\label{eqArmBd4}
	\P_{p_c}\big(|\cluster_G|> \cA^{-4/3}9^k\big)
	\le \P_{p_c}\big(|\cluster_{\Ld}(0)|> \cA^{-4/3}9^k\big)
	\le {C_2}{\cA^{2/3}3^{-k}}. 
	\end{equation}
For the former, on the other hand, we claim that
	\begin{equation}
	\label{eqArmBd1}
	\P_{p_c}\big(H(3^k;G), |\cluster_G|\le \cA^{-4/3}9^k\big)
	\le \cA^{-4/3}3^{k+1}\big(\Gamma(3^{k-1})\big)^2.
	\end{equation}
Indeed, if $|\cluster_G|\le \cA^{-4/3} 9^k$, then there exists $j\in[\frac13 3^k,\frac23 3^k]$ such that $|\partial B_{\cluster_G}(j)|\le \cA^{-4/3} 3^{k+1}$. Denote the first such level by $j$.  
Then, on the right hand side, we get a factor $\Gamma(j)$ (which is bounded by $\Gamma(3^{k-1})$) from the probability of a connection from the origin to level $j$, and $\cA^{-4/3} 3^{k+1}$ times the probability to go from level $j$ to level $3^k$ (each of these probabilities is again bounded above by $\Gamma(3^{k-1})$), which shows \eqref{eqArmBd1}. 

We combine \eqref{eqArmBd3}, \eqref{eqArmBd4}, \eqref{eqArmBd1} with the induction hypothesis, and finally \eqref{eqArmBd2}, to obtain 
	\[	
	\Gamma(3^k) 
	\le \cA^{-4/3}3^{k+1}\left(\frac{\cA}{3^{k-1}}\right)^2+\frac{C_2 \cA^{2/3}}{3^{k}}
	= \frac{3^3\cA^{2/3}+C_2 \cA^{2/3}}{3^k}
	\le \frac \cA{3^k},
	\]
thus proving \eqref{eq Gamma A}. 
\qed

\medskip
While the previous estimates all concern critical percolation, we now turn towards the IIC-measure; and our tool to transfer the results is the construction \eqref{eqPiic}. 

\begin{lemma}[{\cite{KozmaNachm09}}]\label{lemBquant}
Consider percolation in dimension $d>10$. 
There exist $C>0$ such that for all $n\in\N$, $\lambda>1$, 
\[ \P_{\rm\sss IIC}\Big(\frac1\lambda n^2 \le |B(n)| \le \lambda n^2\Big) \ge 1 - \frac C\lambda.\]
\end{lemma}

\proof[Proof of the upper bound.]
We aim to show that $\P_{\rm\sss IIC}\big(|B(n)|>\lambda n^2 \big) \le C\lambda^{-1}$ for all $\lambda>0$, $n\in\N$. 
If $d_{\cluster}(\orig,z)\le n$ and $\orig\conn x$ (for $x,z\in\Zd$), then there exists a vertex $y\in\Zd$ such that 
\[ \{d_{\cluster}(\orig,y)\le n\} \circ \{d_{\cluster}(y,z)\le n\} \circ \{y\conn x\}.\]
By the BK-inequality \eqref{eqBK}, we can bound this from above as follows: 
\begin{align}
	\E_{p_c}[|B(n)|\,\mathbbm 1_{\{\orig\conn x\}}]\,
	&=\sum_{z}\P_{p_c}( d_{\cluster}(\orig,z)\le n, \orig\conn x)\\
	&\le\sum_{y,z}\P_{p_c}\big( \{d_{\cluster}(\orig,y)\le n\} \circ \{d_{\cluster}(y,z)\le n\} \circ \{y\conn x\}\big)\\ 
	&\le\sum_{y,z}\P_{p_c}(d_{\cluster}(\orig,y)\le n)\;\P_{p_c}(d_{\cluster}(y,z)\le n)\,\P_{p_c}(y\conn x). 
\end{align}
Therefore, we get a bound for the conditional probability
\begin{align}
	\E_{p_c}[|B(n)|\,\mid\orig\conn x]\,
	&\le\sum_{x,z}\P_{p_c}(d_{\cluster}(\orig,y)\le n)\;\P_{p_c}(d_{\cluster}(y,z)\le n)\,\frac{\P_{p_c}(y\conn x) }{\P_{p_c}(0\conn x) }. 
\end{align}
The asymptotics \eqref{eqtpf} implies that there is a constant $C'$ such that for all $x$ with $|x-y|\le2|x|$, the ratio $\P_{p_c}(y\conn x) /{\P_{p_c}(0\conn x) }\le C'$, thus 
\begin{align}
	\E_{p_c}[|B(n)|\,\mid\orig\conn x]
	\le C \sum_{x,z}\P_{p_c}(d_{\cluster}(\orig,y)\le n)\;\P_{p_c}(d_{\cluster}(y,z)\le n). 
\end{align}
Finally, we use the upper bound in Lemma \ref{lemBalls} twice to get 
\begin{align}
	\E_{p_c}[|B(n)|\,\mid\orig\conn x]
	\le C' (C_3n)^2.
\end{align}
The finishing touch is provided by Markov's inequality:
\begin{align}
		\P_{p_c}(|B(n)|\ge\lambda n^2\,\mid\orig\conn x)
		\le \frac{C' C_3^2 n^2}{\lambda\,n^2}=C'C_3^2\lambda^{-1}.
\end{align}
Letting $|x|\to\infty$ yields the claim (as $\{|B(n)|\ge\lambda n^2\}$ is a cylinder event). 
\qed

\proof[Proof of the lower bound.] 
For the lower bound, we prove that $\P_{\rm\sss IIC}\big(|B(n)| < \varepsilon n^2 \big) \le C\varepsilon$ for all $\varepsilon=\lambda^{-1}>0$, $n\in\N$. 

If $|B(n)|< \varepsilon n^2$, then there exists some radius $j\in\{\lceil n/2\rceil,\dots,n\}$ such that $|\partial B(0,j)|\le 2\varepsilon n$, and we fix the smallest such $j$. 
Then we condition on $\{B(j)=A\}$ for any ``$j$-admissible'' subgraph $A$, which is any finite subgraph $A$ of $\Ld$ containing $\orig$ s.\ t.\
\begin{itemize}
\item $\P_{p_c}(B(j)=A)>0$,
\item $|\partial A|\le 2\varepsilon n$, where $|\partial A|$ denote the number of vertices at maximal graphical distance from $\orig$ 
\item $|\{y\colon d_A(\orig,y)=k\}|> 2\varepsilon n$ for $k=\lceil n/2\rceil,\dots,j-1$ (to make sure that $j$ is the ``first'' level satisfying the above property).
\end{itemize} 
This yields 
\begin{align}
	\P_{p_c}\big(|B(n)|<\varepsilon n^2,0\conn x \big) 
	\,&\le \sum_{j=n/2}^n \sum_{A}\P_{p_c}\big(B(j)=A,0\conn x \big) \nonumber \\
	&= \sum_{j=n/2}^n \sum_{A}\P_{p_c}\big(0\conn x\mid B(j)=A \big)
		\;\P_{p_c}\big(B(j)=A\big),  \label{eqBjnBd}
\end{align}
where the sum is over all $j$-admissable $A$. 
For any such $A$, we get
\[ \P_{p_c}(0\conn x\mid B(j)=A)
	\le \sum_{y\in\partial A} \P_{p_c}(y\conn x \text{ with a path avoiding $A\setminus\partial A$ }\mid B(j)=A). \]
However, since $\{y\conn x \text{ with a path avoiding $A\setminus\partial A$}\}$ only depends on the edges with both endpoints outside $A\setminus\partial A$  and $\{B(j)=A\}$ only depends on the edges with both endpoints in $A$, the two events are independent, and 
\begin{align*}
	\P_{p_c}(0\conn x\mid B(j)=A)
	\, &\le \sum_{y\in\partial A} \P_{p_c}(y\conn x \text{ with a path avoiding $A\setminus\partial A$})\\
	&\le \sum_{y\in\partial A} \P_{p_c}(y\conn x)
	\le \sum_{y\in\partial A} C |y-x|^{d-2}, 
\end{align*}
where the last bound uses \eqref{eqtpf}.
Assuming that $x$ is far away from the origin (again $|x-y|\le 2|x|$ suffices), then there is a constant $C'>0$ such that 
\[ \P_{p_c}(0\conn x\mid B(j)=A) \le C' \,|\partial A|\,|x|^{2-d} \le  C' \varepsilon n \,|x|^{2-d}.\]
Furthermore, we have that 
\[ \sum_{j=n/2}^n \sum_{A}\P_{p_c}\big(B(j)=A\big)\le P_{p_c}\big(\partial B(n/2)\neq\varnothing\big).\]
Plugging the previous two bounds in \eqref{eqBjnBd}, we get
\begin{align*} \P_{p_c}\big(|B(n)|<\varepsilon n^2,0\conn x \big) 
	&\le C'\varepsilon n |x|^{2-d}\sum_{j=n/2}^n \sum_A\P_{p_c}(B(j)=A)\\
	&\le C'\varepsilon n |x|^{2-d}\P_{p_c}\big(\partial B(n/2)\neq\varnothing\big), 
\end{align*}
and now we use the upper bound in Lemma \ref{lemArmExponent} to further bound 
\[ \P_{p_c}\big(|B(n)|<\varepsilon n^2,0\conn x \big) \le C''\varepsilon |x|^{2-d}\]
for a constant $C''>0$.  
Finally, letting $|x|\to\infty$ and using \eqref{eqPiic} along with the two-point function estimate \eqref{eqtpf} yields the desired result. 
\qed

\medskip

In order to prove that $\df(\cluster)=2$ for the incipient infinite cluster, we combine the previous lemma with the following general criterion: 
\begin{lemma}[{Lemma 3.2 in \cite{Cames15}}]\label{lemCriterion}
Let $(Z_n)_{n\in\N}$ be a sequence of positive random variables such that $Z_1\le Z_2\le\dots$. Suppose there are constants $\alpha,\mu,C>0$ such that for all $\lambda>0$ and $n\in\N$, we have 
\begin{equation}\label{eqLambdaBd}
\P(\lambda^{-1}n^\alpha \le Z_n \le \lambda n^\alpha)\ge 1-C(\log \lambda)^{-1-\mu}.
\end{equation}
Then 
\[\P\left(\lim_{n\to\infty}\frac{\log Z_n}{\log n}=\alpha\right)=1.\]
\end{lemma}

\proof
We abbreviate $Y_n:=\log Z_n / \log n$, and claim that it is sufficient to prove  
\begin{equation}\label{eqConvergence2n}
	\lim_{k\to\infty}Y_{2^k}=\alpha \qquad \P-a.s.
\end{equation}
Indeed, for $n\in\N$, we choose $k=k(n)\in\N$ such that $2^{k-1}\le n\le 2^k$, and use the monotonicity of the sequence $(Z_n)_n\in\N$ to bound 
\[ Y_{2^{k-1}}\frac{k-1}{k}= \frac{\log Z_{2^{k-1}}}{\log 2^{k}}\le\frac{\log Z_n}{\log n}
	\le\frac{\log Z_{2^{k}}}{\log 2^{k-1}}=Y_{2^{k}}\frac{k}{k-1}, 
\]
and then use \eqref{eqConvergence2n} to conclude the claim. 

In order to prove \eqref{eqConvergence2n}, we define 
\begin{equation*}
	\varepsilon_k:= k^{\frac{1+\mu/2}{1+\mu}-1},
	\qquad \lambda_k:= 2^{k\varepsilon_k},
\end{equation*}
and note that $\varepsilon_k>0$, $\lambda_k>1$ for all $k\ge1$, and $\lim_{k\to\infty}\varepsilon_k=0$. 
Then, using \eqref{eqLambdaBd}, 
\begin{align*}
	\sum_{k=1}^\infty\P\big(|Y_{2^k}-\alpha|> \varepsilon_k\big)
	\,&= \sum_{k=1}^\infty\P\big(|\log Z_{2^k}- \log(2^{k\alpha})|> \log\lambda_k\big)\\
	\,& = \sum_{k=1}^\infty\P\big(Z_{2^k}< \lambda_k^{-1} 2^{k\alpha}\big)+\P\big(Z_{2^k}> \lambda_k 2^{k\alpha}\big)\\
	& \le C\sum_{k=1}^\infty \frac1{(\log\lambda_k)^{1+\mu}}
	 = \frac C{(\log2)^{1+\mu}}\sum_{k=1}^\infty \frac1{k^{1+\mu}}<\infty. 
\end{align*}
Hence, the Borel-Cantelli lemma implies that 
\[ \P\big(|Y_{2^k}-\alpha|> \varepsilon_k \text{ for infinitely many }k\big)=0, \]
which proves \eqref{eqConvergence2n}.
\qed

\proof[Proof of $\df(\cluster)=2$.] 
We apply Lemma \ref{lemCriterion} with $\P$ being the IIC-measure $\P_{\rm\sss IIC}$, $\alpha=2$, $Z_n=B(n)$ and apply Lemma \ref{lemBquant} to get the desired result. 
\qed

\subsection{The spectral dimension}
Control of the return probability of random walk needs two ingredients. The first one is control of the the volume growth, which is achieved in Lemma \ref{lemBquant}. The second ingredient is control of the \emph{effective resistance}. The connection between these two ingredients and random walk behaviour is in the folklore of studying random walks, see in particular Kumagai and Misumi \cite{KumagMisum08} for results in our context. Kozma and Nachmias \cite{KozmaNachm09} prove a quantitative estimate on the lower bound on the effective resistance between $\orig$ and $\partial B(n)$, and then apply a readily tailored theorem of Barlow at al.\ \cite{BarloJaraiKumagSlade08} to deduce that $\ds=4/3$. Another consequence of this theorem is that the escape time exponent equals $\beta=3$ $\Piic-$almost surely (precisely as for the IIC on trees).

\subsection{The mass dimension}
Already van der Hofstad and J\'arai \cite{HofstJarai04} showed that 
\[ {\mathbb{E}_{\scriptscriptstyle \IIC}}|\cluster\cap Q(n)|\approx n^4.\] 
From this, we can prove that $\df(\cluster)\le 4$ rather straightforwardly via Markov's inequality. The challenge is to prove a complementing lower bound, which was achieved by Cames van Batenburg \cite{Cames15} using quantitative bounds on the extrinsic one-arm exponent \cite{KozmaNachm11}. 

Mind that the escape time exponent as defined in \eqref{eqDefBeta} determines the rate at which a random walk leaves a ball of intrinsic distance $n$. Unlike on $\Zd$, the extrinsic and intrinsic distances are not equivalent on the IIC-cluster, and we therefore consider a modified critical exponent $\beta'$ as 
\begin{equation}
	E\big[\inf\{n\in\N\colon S_n\in\partial Q(n)\}\big]=n^{\beta'+o(1)}.
\end{equation}
With van der Hofstad and Hulshof \cite{HeydeHofstHulsh14a} we proved that $\beta'=6$ for $\Piic-$almost all realizations $\omega$. This should be contrasted against $\beta=3$ explained before. This means that the random walk needs order $n^3$ steps to leave the intrinsic ball $B(n)$, but it needs $n^6$ steps in order to leave $Q(n)$. 
The factor 2 between these two exponents is not a coincidence: in high dimensions, the spatial dependency between different parts of a critical cluster is rather weak; in fact so weak that geodesic paths (w.r.t.\ graph distance) are embedded into $\Zd$ similar to a random walk path, and thus the graph distance between $\orig$ and $\partial Q(n)$ is of the order $n^2$. 

\section{Discussion and outlook}
A number of pressing challenges were mentioned en passant, most notably the identification of dimensions of critical percolation clusters in lower dimension. However, in the following we want bring forward two lines of further research that might be within reach with current techniques. 

\begin{enumerate}
\item Identify discrete Hausdorff dimension and packing dimension of critical percolation clusters in high dimension. 
Also for other ``natural'' random subsets of $\Zd$. So far only results by Barlow and Taylor \cite{BarloTaylo92} and Georgiou et al.\ \cite{GeorgKhoshKimRamos18} for the range of (generalised) random walks.

\item Known cases of discrete dimension all deal with subsets of $\Zd$, and also the focus of the present account is on subsets of the hypercubic lattice. However, there is no obvious need to stick to the lattice setup here---fractal and spectral dimension are meaningful for any locally-finite connected graph, and the others require an embedding of the vertices in some metric space, and $\Zd$ might appear as an unnecessary limitation. 

From a geometric point of view, it might be more natural to focus on discrete subsets of $\Rd$. Instead of lattice percolation, one might investigate the geometric properties of (critical) continuum percolation clusters. A suitable candidate is the \emph{random connection model}, where vertices are given as a Poisson point process in $\Rd$, and two vertices are linked by an edge with probability depending on the Euclidean distance between the vertices. The critical behaviour of the random connection model in high dimensions has recently been identified \cite{HeydeHofstLastMatzk19}, paving the way to an investigation of the continuum incipient infinite cluster and its dimension(s). 
\end{enumerate}

\subsection*{Acknowledgement.} The author thanks Martin Barlow and Steffen Winter for providing references and for comments on an earlier version of the manuscript. 

%
\def\cprime{$'$}
\providecommand{\bysame}{\leavevmode\hbox to3em{\hrulefill}\thinspace}
\providecommand{\MR}{\relax\ifhmode\unskip\space\fi MR }
\providecommand{\MRhref}[2]{%
  \href{http://www.ams.org/mathscinet-getitem?mr=#1}{#2}
}
\providecommand{\href}[2]{#2}

\end{document}